\documentclass[reqno]{amsart}
\usepackage{amsfonts,amsmath}
\usepackage{amssymb}
\usepackage[cp1251]{inputenc}
\usepackage{graphicx}
\newtheorem{theorem}{Theorem}%[section]

 \theoremstyle{definition}
 \newtheorem{defn}{Definition}[section]
\UseRawInputEncoding
\usepackage[russian,english]{babel}
\numberwithin{equation}{section}
\def\dfrac#1#2{\displaystyle{#1\over #2}}

\def\bV{{\bf v}}

\def\n{\rho}

\def\Div{\mbox{div}\,}
\def\Rot{\mbox{rot}\,}

\def\bB{{\bf B}}

\def\bE{{\bf E}}
\begin{document}

\title[Internal free boundary problem] {Internal free boundary problem for cold plasma equations}

\author[Gargyants]{L.V. Gargyants$^1$}
\author[Konovalova]{A.V. Konovalova$^2$}
\author[Rozanova]{ O.S. Rozanova$^3$}

\address[1]{Bauman Moscow State Technical University}%{City, Country}{email}
\address[2,3]{Lomonosov Moscow State University}
%\email{rozanova@mech.math.msu.su}

\begin{abstract} For the system of cold plasma equations describing the motion of electrons in the field of stationary ions, we consider the Riemann problem posed at an impenetrable interface between two media. These media differ in the magnitude of the constant ion field. The interface between the media is assumed to be free. Its position is determined from the generalized Rankine-Hugoniot conditions and the stability condition, that is, the intersection of Lagrangian particle trajectories at the interface.
\end{abstract}
\maketitle

%\begin{history}
%\received{(Day Month Year)}
%\revised{(Day Month Year)}
%\accepted{(Day Month Year)}
%\comby{(xxxxxxxxxx)}
%\end{history}

\section{Introduction}

 The system of hydrodynamic of electron liquid, together with Maxwell's equations, has the following form:
\begin{equation*}
\label{base1}
\begin{array}{l}
 \n_t + \Div(\n \bV)=0\,,\quad
\bV_ t + \left( \bV \cdot \nabla \right) \bV
=\dfrac {e}{m} \, \left( \bE + \dfrac{1}{c} \left[\bV \times  \bB\right]\right),\vspace{0.5em}\\
\dfrac1{c}  \bE_t = - \dfrac{4 \pi}{c} e \n \bV
 + {\rm rot}\, \bB\,,\quad
\dfrac1{c}  \bB_ t  =
 - {\rm rot}\, \bE\,, \quad \Div \bB=0\,,
\end{array}
\end{equation*}
where $e, m$ are the charge and mass of the electron (here the electron charge has a negative sign: $ e <0 $),
$ c $ is the speed of light;
$ \n,  \bV $ are the density and the velocity vector of
electrons;
$ \bE, \bB $ are the vectors of electric and magnetic fields, $x\in{\mathbb R}^3,$ $t\ge 0$, $\nabla$, $\rm div$, $\Rot$ are the gradient, divergence and vorticity with respect to the spatial variables.
The system is often called the equations of hydrodynamics of "cold" plasma (see, for example,  \cite {ABR}, \cite {david72}).

Interest in cold plasma models has always been high, as they represent a reduction and therefore a possible simplification of the equations describing plasma in general. However, this interest has recently increased further due to the practical feasibility of creating accelerators on the wake wave (a review can be found in \cite{Shep13}, \cite{Malka}, see also \cite{GV}, \cite{KOS}), particularly those used in medical device (\cite{MKS} and references therein).

Cold plasma (or electron liquid) of  is a highly unstable medium, in which even small perturbations of a constant state lead to a loss of smoothness in the solution.
Such phenomena in the physically natural three-dimensional case are practically impossible to study analytically, and even numerical analysis is associated with significant difficulties (see, for example, \cite{Ch_book} and references therein).

However, if we consider the so-called plane plasma oscillations, in which each component of the solution depends on only one spatial variable, the system is significantly simplified and analytical results are possible.
Namely,in this in the one-dimensional in space case
$\bV=(V,0,0)$, $\bE=(E,0,0)$, $\bB\equiv 0$, e.g.\cite{Ch_book}. In dimensionless form it can be written as
 \begin{equation}
\begin{array}{c}
\n_t +
\left(\n\, V \right)_x
=0,\quad
V_t +
V  V_x =  - E, \quad
E_t = \n\, V.
\end{array}
\label{Cold}
\end{equation}
For smooth solutions  (\ref{Cold}) implies
%\begin{equation}
$
\left(\n +
 E_x \right)_t = 0,$
%\end{equation}
 it follows
\begin{equation}\label{n}
\n=N(x)-E_x.
\end{equation}
The function $N(x)\ge 0$ is called the background density of electrons.

If we substitute \eqref{n} to the last equation of  \eqref{Cold}, we can exclude the density and obtain system
\begin{equation}
 \begin{array}{c}
E_t + V  E_x
=N(x) V,\quad
V_t +
V  V_x =  - E,
\end{array}
\label{EV}
\end{equation}
very convenient for the analysis of smooth solution. Indeed, system is non-strictly hyperbolic and all dynamics can be studied along one Lagrangian characteristics \cite{Nonstrictly}. It $E$ is known, then $\n$ can be found from \eqref{n}.

The following is known about system \eqref{EV}.

1. Let $(E_0(x), V_0(x))$ be $C^1$ initial data for \eqref{EV} and $N(x)=c={\rm const}>0$. Then there exist nontrivial initial data such that the solution to the Cauchy problem keeps smoothness for all $t>0$. Specifically,
the class of these initial data is determined by condition \cite{RCh_ZAMP}
\begin{equation*}
(V_0')^2+ 2 E_0'-c<0.
\end{equation*}

2. If $N(x)>0$ is not equal to a constant, then any nontrivial solution blows up in a finite time \cite{PhD_doping}.

This motivates to consider the initial data that are initially discontinuous and study the evolution of the discontinuity at least for simplest situations.  In \cite{Riemann_Roz} the Riemann initial data and constant background density were considered. Note that constructing the shock wave requires a conservative form of the system.
It was established that

1. the rarefaction wave and shock wave alternate periodically;

2. The shock wave is strongly singular (i.e., it contains the density in the form of  delta function).

3. The rarefaction wave is generally not unique, and uniqueness requires the use of unconventional techniques, such as requiring a minimum energy within the rarefaction wave.

In this paper, we continue to study the Riemann problem for a one-dimensional cold plasma in the case where the background density initially exhibits a strong discontinuity. In other words, we want to study the dynamics of the discontinuity between two media, which can be considered an internal free boundary. Therefore, the problem can be considered a variant of the Riemann problem.

Specifically, we consider \eqref{Cold}
in  $\Omega_-(t): x<\Phi(t)$ or   $\Omega_+(t): x>\Phi(t)$,
separated by an impenetrable boundary $x=\Phi(t)$, $\Phi(0)=0$, the position of which is unknown in advance. For each of these regions, we assign a different constant value of the background density:
\begin{equation*}\label{nn}
N(x)=n_\pm={\rm const}>0.
\end{equation*}

As follows from \eqref{n},
\begin{equation}
 \n = n_\pm -E_x.
%\dfrac{\partial  E }{\partial x}.
%\label{3gl4}
\label{Kn}
\end{equation}

Here and below for all functions $f$, having one-sided limits $f_\pm$ at a point $x_0\in\mathbb R$ we denote by $[f]=f_+-f_-$ the value of jump at $x_0$.

For  \eqref{EV} we prescribe initial conditions
\begin{equation}\label{Data}
(V,E)|_{t=0}=(V^0_-+[V]^0 \Theta (x),\,\,  E^0_-+[E]^0 \Theta (x)),\quad x \in \mathbb R,
\end{equation}
and \eqref{Kn} implies
\begin{equation*}\label{K3}
\n\big|_{t=0}=n_\pm- [E]^0\delta(x).
\end{equation*}
Here \begin{equation*}\label{K3}
V^0_-, \,[V]^0, \, E^0_-, \, [E]^0<0,
\end{equation*}
are constant, $[.]^0$ is a jump at zero.

The problem consists of finding a solution to the problem \eqref{EV}, \eqref{Data} in each of the regions $\Omega_\mp(t)$, as well as finding the boundary
$x=\Phi(t)$ between these regions based on the  admissibility condition and the Rankine-Hugoniot conditions.

There are two possible situations when $t\in (0,T)$, $T>0$ is small enough.

\medskip

1. $V^0_->V^0_+$, when the free boundary can be interpreted as a line of singular strong discontinuity. The characteristics (Lagrangian trajectories) intersect at $x=\Phi(t)$, and mass accumulates at the discontinuity according to the generalized Rankine-Hugoniot conditions. The admissibility condition for a singular shock wave coincides with the geometric entropy condition:
\begin{eqnarray}
\label{accept}
 \min\{V_-, V_+\} \le  \dot \Phi(t) \le \max\{V_-, V_+\},
\end{eqnarray}
meaning that characteristics from both sides come to the shock.

\medskip

2. $V^0_-<V^0_+$, when a rarefaction wave forms to the right and left of $x=\Phi(t)$. It is separated from the state independent of $x$ by the curves $x=x_-(t)$ and $x=x_+(t)$, respectively. Thus, the $\Omega_-$ region consists of $\Omega^1_-:x<x_-(t)$ and
$\Omega^2_-: x_-(t)<x<\Phi(t)$, while $\Omega_+$ consists of $\Omega^1_+:x> x_+(t)$ and
$\Omega^2_+: x_+(t)>x>\Phi(t)$.

\medskip
The characteristics $x_-(t)$ and $x_+(t)$, outgoing from the  point $x=0$ have an infinite number of intersection points. Thus, rarefaction regions and shock waves alternate sequentially, and their initial and final points can be found from a transcendent equation. Below, we construct a singular shock wave and a rarefaction wave in each interval. Note that they are found based on different forms of the system.

Note that the construction of a singular wave at the interface between two media is practically identical to the situation of a constant background density described in \cite{Riemann_Roz}, whereas instead of a single rarefaction wave, a structure is formed consisting of either two rarefaction waves separated by a singular shock wave, or one rarefaction wave and a singular shock wave.

\section{Lagrangian trajectories (characteristics)}\label{S2}

The characteristics of the system from both sides of the interface of media can be found from \eqref{EV}, which reads as
\begin{equation}\label{K1}
V_t+V V_x=-E, \quad E_t+V E_x=n_\pm V
\end{equation}
in every $\Omega_\pm$.
The characteristics found as a solution to the system
 \begin{equation*}\label{Kchar}
 \dfrac {d V}{dt}= -E,\quad \dfrac {d E}{dt}=n_\pm  V,\quad \dfrac {d x}{dt}= V,
 \end{equation*}
 that is for the Lagrangian trajectories, starting from the point $x_0$,
\begin{eqnarray}\label{KV}
   &V_\pm(t)=-\dfrac{E^0_\pm}{\sqrt{n_\pm}} \sin \sqrt{n_\pm} t + V^0_\pm \cos \sqrt{n_\pm} t, \\
   &  E_\pm(t)=V^0_\pm\sqrt{n_\pm} \sin \sqrt{n_\pm} t + E^0_\pm \cos \sqrt{n_\pm} t,\label{KE}\\
   &x_\pm(t)=\dfrac{V^0_\pm}{\sqrt{n_\pm}} \sin \sqrt{n_\pm} t + \dfrac{E^0_\pm}{n_\pm} ( \cos \sqrt{n_\pm} t -1)+ x_0. \label{Kx}%\quad x_0=0.
 \end{eqnarray}

Let us denote $T_*>0$ as the smallest positive moment of time when the characteristics $x_\pm(t)$ intersect, i.e. $x_-(t)=x_+(t).$ At this moment of time, the shock wave transforms into a rarefaction wave or vice versa, depending on the Riemann data.

Note that if $\sin{\sqrt{n_\pm}T_*}=0,$ and $\cos{\sqrt{n_\pm}T_*}\ne 1,$ then $x_-(T_*)=x_+(T_*)$ if and only if $\frac{E_+^0}{n_+}=\frac{E_-^0}{n_-}.$

\section{Construction of a singular shock wave}\label{S3}

Assume that  $V^0_->V^0_+$, that is, the solution to the Riemann problem starts from a shock wave.

To construct a  shock wave, the system must be written in divergent form:
\begin{equation}\label{K4}
\n_t + (V \n)_x=0, \qquad
\left(\frac{\n V^2}{2}+\frac{E^2}{2}\right)_t +\left(\frac{\n V^3}{2}\right)_x =0,
\end{equation}
corresponding to the laws of conservation of mass and total energy
(for example, \cite {FrCh}).

Since an initial discontinuity in the component $E$
implies the delta-singularity of $\n$ (see \eqref{Kn}), then we have to use a concept of the singular shock \cite{Shelkv} and seek for the solution as
\begin{eqnarray}\label{K19}
   V(t,x)&=&V_-(t,x)+[V(t,x)]|_{x={\Phi}(t)} \Theta(x-{\Phi}(t)), \\  E(t,x)&=&E_-(t,x)+[E(t,x)]|_{x={\Phi}(t)} \Theta(x-{\Phi}(t)),\label{K20}
   \\ \rho(t,x)&=& \hat\rho(t,x)%+[\rho(t,x)]|_{x={\Phi}(t)} \Theta(x-{\Phi}(t))
   +e(t)\delta(x-{\Phi}(t)),\label{K21}
    \end{eqnarray}
    where $[f]=f_+-f_-$, $f_\pm $ are differentiable functions with one-sided limits, $t\ge 0$, $x\in \mathbb R$, $ \hat \rho(t,x)=n_\pm-\{E_x(t,x)\}$, $\{E_x\}$ is the derivative of $E$ at those points where it exists in the usual sense, $e(t):=e(t,{\Phi}(t))$, $e(t)=-[E(t,x)]|_{x={\Phi}(t)}$.

    \subsection{Definition of a generalized strongly singular solution}

 Starting from the divergent form \eqref {K4}, we define a generalized strongly singular solution according to  \cite{Shelkv}, \cite{NRS}.

\begin{defn}
The triple of distributions $ (V, E, \n) $, given as \eqref {K19} - \eqref {K21} and the curve $ \gamma $, given as $ x = \Phi (t), $ $ \Phi (0) = 0 $, $ \Phi (t) \in C ^ 1 $, is called a generalized singular solution of the problem \eqref{K4},
\begin{eqnarray*}\label{K30}
(V,E,\n)|_{t=0}=\\(V^0_-(x)+[V(x)]^0 \Theta (x),\, E^0_-(x) +[E(x)]^0 \Theta(x),\, \n^0(x)=\hat \n^0(x)+e^0\delta(x)),\nonumber
\end{eqnarray*}
if for all test functions $\phi(t,x)\in \mathcal{D}({\mathbb R}\times [0,\infty))$
\begin{eqnarray*}\nonumber
\int\limits_0^\infty\int\limits_{\mathbb R} \hat \n  (\phi_t+V \phi_x) dx dt +\int\limits_{\gamma} e(t) \frac{\delta \phi (t,x)}{\delta t} \frac{dl}{\sqrt{1+(\dot\Phi(t))^2}}+\\
\int\limits_{\mathbb R} \hat \n^0(x) \phi(0,x) dx  + e(0) \phi(0,0)  =0,\label{K190}\\
\int\limits_0^\infty\int\limits_{\mathbb R} \left(( \frac{\hat \n V^2}{2}+E^2)\phi_t+ \frac{\hat \n V^3}{2} \phi_x\right) dx dt +\int\limits_{\gamma} \frac{e(t)(\dot \Phi(t))^2}{2} \frac{\delta \phi (t,x)}{\delta t}  \frac{dl}{\sqrt{1+(\dot\Phi(t))^2}}+\nonumber\\
\int\limits_{\mathbb R} \left(\frac{\hat \n^0(x) (V^0(x))^2}{2} + (E^0(x))^2 \right)\phi(0,x) dx  + \frac{e(0)(\dot \Phi(0))^2}{2} \phi(0,0) =0,\label{K200}
\end{eqnarray*}
 where $\int\limits_{\gamma} \cdot dl$ is the curvilinear integral along the curve $ \gamma $, the delta-derivative $\frac{\delta \phi (t,x)}{\delta t}\big|_{\gamma} $ is defined as the tangential derivative on the curve $ \gamma $, namely
$$
\frac{\delta \phi (t,x)}{\delta t}\big|_{\gamma} =\left(\frac{\partial \phi (t,x)}{\partial t}+ \dot\Phi(t) \frac{\partial \phi (t,x)}{\partial x}\big|_{\gamma} \right)\big|_{\gamma}=\frac{d \phi (t,\Phi(t))}{d t}= \sqrt{1+(\dot\Phi(t))^2} \frac{\partial \phi (t,x)}{d {\bf l}},
$$
where ${\bf l}=(-\nu_2, \nu_1)=\frac{(1, \dot\Phi(t))}{\sqrt{1+(\dot\Phi(t))^2}}$ is a unit vector tangent to $\gamma$.

\end{defn}
The action of the delta function $ \delta (\gamma) $ concentrated on the curve $ \gamma $ on the test function is defined according to \cite {Kanwal},
as $$
(\delta(\gamma),\phi(t,x))=\int\limits_{\gamma} \phi(t,x) \frac{dl}{\sqrt{1+(\dot\Phi(t))^2}},
$$
where $\phi(t, x)\in \mathcal{D}({\mathbb R}\times [0,\infty))$.

\bigskip

    Direct adaptation of Rankine-Hugoniot conditions to the case of a singular shock   (\cite{Shelkv}, \cite{Riemann_Roz}) gives the following theorem.
\begin{theorem} Let the domain $\Omega\subset {\mathbb R}^2$ be divided by a smooth curve $\gamma=\{(t,x): x=\Phi(t) \}$ into left and right parts
$\Omega_\mp$. Let the triple of distributions $(V,E,n)$ \eqref{K19} -- \eqref{K21} and the curve $\gamma$ be a strongly singular generalized solution for the system \eqref{K4}. Then this solution satisfies the following analogue of the Rankine-Hugoniot conditions
 \begin{eqnarray}\label{RH1}
    e'(t)&=&\left(-[ \hat \rho V]+[\hat \rho] \Phi'(t)\right)\big|_{x=\Phi(t)},\\
   \label{RH2}
  \frac{d}{dt}\frac{ e(t) ( \Phi'(t))^2}{2}&=&\left(-\left[\frac{\hat \rho V^3}{2}\right]+\left[\frac{\hat \rho V^2+E^2}{2}\right] \Phi'(t)\right)\big|_{x=\Phi(t)}.
 \end{eqnarray}
 \end{theorem}

 \bigskip

 {\it The proof} of the theorem is completely analogous to the proof of a similar result for the case $n_-=n_+$, contained in \cite{Riemann_Roz}, so we do not present it.

 The theorem allows to find the position of the interface in dependence on the values of the components from both sides of the interface.

\subsection{Case of the  Riemann data}
The Riemann data \eqref{Data} do not depend of  $x$, therefore the solution  to  \eqref{EV}, \eqref{Data} in $\Omega_\pm$, do not depend  of the spatial variables, as well, and are given by \eqref{KV}, \eqref{KE}.

\begin{eqnarray*}\label{K19a}
   V(t,x)&=&V_-(t)+[V(t)]|_{x={\Phi}(t)} \Theta(x-{\Phi}(t)), \\  E(t,x)&=&E_-(t)+[E(t)]|_{x={\Phi}(t)} \Theta(x-{\Phi}(t)),\label{K20a}
   \\ \rho(t,x)&=& n_-+(n_+-n_-) \Theta(x-{\Phi}(t))+e(t)\delta(x-{\Phi}(t)),\label{K21a}
    \end{eqnarray*}
Assume that the shock starts at the point $t=T_0\ge 0$.

Therefore, ~\eqref{RH1}  takes the form
%\begin{eqnarray}
%&& e'(t)=-(-\sqrt{n_+}E_+^0\sin{\sqrt{n_+}t}+n_+V_+^0\cos{\sqrt{n_+}t})+(-\sqrt{n_-}E_-^0\sin{\sqrt{n_-}t}+\nonumber\\
%&&+n_-V_-^0\cos{\sqrt{n_-}t})+(n_+-n_-)\Phi'(t),\label{esw}
%\end{eqnarray}
\begin{eqnarray}
&& e'(t)=-[nV(t)]\big|_{x=\Phi(t)}+[n]\big|_{x=\Phi(t)}\Phi'(t),\label{esw}
\end{eqnarray}
from where, taking into account the initial conditions $e(T_0)=-[E(T_0)],\;\Phi(T_0)=\Phi_0,$ we obtain
\begin{equation}\label{ee}
  % e(t)=-[E(t)]-[E]^0 +[n](\Phi(t)-\Phi_0)\ge 0.
   e(t)=-[E(t)]\big|_{x=\Phi(t)}+[n]\big|_{x=\Phi(t)}(\Phi(t)-\Phi_0).
\end{equation}

Note that \eqref{KV}, \eqref{KE} implies
\[
[\hat\n V^2+E^2]\big|_{x=\Phi(t)}=[\hat\n V^2+E^2]^0=K=\rm const,
\]
where $\hat\n= \rho_-(t,x)+[\rho(t,x)]|_{x={\Phi}(t)} \Theta(x-{\Phi}(t))$,  the regular component of density.

Therefore, ~\eqref{RH2} can be written as
\begin{equation}\label{ePhi}
    ({ e(t) (\Phi'(t))^2})'=-\left[{\hat \rho V^3(t)}\right]\big|_{x=\Phi(t)}+{K}\Phi'(t).
\end{equation}
Substitute \eqref{ee} to \eqref{ePhi} and obtain a nonlinear second order equation with respect to $\Phi(t)$ ,
\begin{eqnarray}
\label{Phidd}
&&\Phi'(t) \Big( 2 (\alpha \Phi(t)+\beta(t))  \Phi''(t) + \alpha(\Phi'(t))^2 + \beta'(t) \Phi'(t)-K  \Big)=- \sigma(t),
\end{eqnarray}
with $\alpha=n_+-n_-$, $\beta(t)=-[E(t)]\big|_{x=\Phi(t)}-\alpha \Phi_0$, $\sigma(t)=\left[{\hat \rho V^3}\right]\big|_{x=\Phi(t)}$,
which should be supplemented with two boundary conditions
\begin{equation*}\label{bound}
\Phi(T_0)=\Phi_0,\quad \Phi(T_*)=x_+(T_*)=x_-(T_*).
\end{equation*}
We see that the equation for determining $\Phi(t)$ is much more complicated than in the case of a constant background density, when $\alpha=0$.

\medskip
\section{ Construction of a rarefaction wave}\label{S4}

The solution to the  problem \eqref{EV}, \eqref{Data} in $\Omega^1_\pm$, independent of the spatial variables, as in the case of shock wave, has the form~\eqref{KV}, \eqref{KE}. However, the domains $\Omega^2_\pm$ are not occupied by the characteristics, and
we need to construct solution, that is continuous at least on each side of the interface $x=\Phi(t)$. For the case of constant background density, it is possible to construct a continuous solution everywhere between diverging characteristics $x_-(t)$, $x_+(t)$ but it can be shown that this cannot be achieved in the case of contact between media with different background densities. Therefore, the domains $\Omega_\pm$  in which a continuous solution is constructed must be connected by a discontinuity, following the principles of constructing a singular shock wave.

Thus, in the  domains $\Omega^2_\pm$ , we must construct a solution to system \eqref{K1} and find the location of the interface $x=\Phi(t)$ based on the generalized Rankine-Hugoniot conditions.

Essentially, in each domain $\Omega^2_\pm$, the solution is uniquely found as a function affine in the spatial variable, and the problem consists of choosing the correct interface.

So,
\[
(V,E)=
\begin{cases}
&(V_-(t),E_-(t)),\quad x<x_-(t),\\
&(V_\pm^r,E_\pm^r)=(a_\pm(t)x+b_\pm(t), c_\pm(t)+d_\pm(t)),\quad x \in \Omega^2_\pm,\\
&(V_+(t),E_+(t)),\quad x>x_+(t)
\end{cases}
\]
Assume that the rarefaction wave starts at $t=T_0=0$ and  it is replaced by a shock  at $t=T_*$. Otherwise, we change the time to $t-T_0$.

The coefficients $a_\pm(t), c_\pm(t)$ satisfy the system of differential equations:
\[
\dot a_\pm=-a_\pm^2-c_\pm, \quad \dot c_\pm= a_\pm (n_\pm- c_\pm).
\]
We solve this system subject to condition $a_\pm(t), c_\pm(t) \to \infty$ at $t\to 0$ and at $t\to T_*$ and obtain
\begin{eqnarray*}
&&a_\pm(t)=\dfrac{\sqrt{n_\pm}(B_\pm\cos{\sqrt{n_\pm}t}-\sin{\sqrt{n_\pm}t})}{B_\pm\sin{\sqrt{n_\pm}t}+\cos{\sqrt{n_\pm}t}-1},\\ &&c_\pm(t)=\dfrac{n_\pm(B_\pm\sin{\sqrt{n_\pm}t}+\cos{\sqrt{n_\pm}t})}{B_\pm\sin{\sqrt{n_\pm}t}+\cos{\sqrt{n_\pm}t}-1},
\end{eqnarray*}
with
\[
B_\pm=\frac{1-\cos{\sqrt{n_\pm}T_*}}{\sin{\sqrt{n_\pm}T_*}}=\tan{\frac{\sqrt{n_\pm}T_*}{2}}.
\]
Note that if $ \sin{\sqrt{n_\pm}T_*}= 0$, then the limit pass shows that $B_\pm=0$.

Further, from the continuity on the boundary of $\Omega_\pm^1$ and $\Omega_\pm^2$ we have
\begin{eqnarray*}
&&b_\pm(t)=V_\pm(t)-a_\pm(t)x_\pm(t)=
\frac{(B_\pm E_\pm^0-\sqrt{n_\pm}V_\pm^0)(\cos{\sqrt{n_\pm}t}-1)}{\sqrt{n_\pm}(B_\pm\sin{\sqrt{n_\pm}t}+\cos{\sqrt{n_\pm}t}-1)}\\
&&d_\pm(t)=E_\pm(t)-c_\pm(t)x_\pm(t)
=\frac{(B_\pm E_\pm^0-\sqrt{n_\pm}V_\pm^0)\sin{\sqrt{n_\pm}t}}{B_\pm\sin{\sqrt{n_\pm}t}+\cos{\sqrt{n_\pm}t}-1}.
\end{eqnarray*}
From ~\eqref{Kn} we obtain $\rho=n_\pm-c_\pm(t).$

Note that for existence of the shock we have to require $e(t)\ge 0$, this implies
\begin{equation}\label{cmcp}
c_-(t) \Phi(t)+ d_-(t) \ge c_+(t) \Phi(t)+ d_+(t), \qquad t\in (0,T_*).
\end{equation}

Then \eqref{RH1}, \eqref{RH2} imply
 \begin{eqnarray}
 &e'\,&=\,[n-c]\,\Phi' -[(n-c)a]\,\Phi-[(n-c)b],\label{eR}
  \\
 & (e (\Phi')^2)'&\,=\,-[(n-c) a^3]\, \Phi^3 + 3 [(n-c) a^2 b]\, \Phi^2 - 3 [(n-c) a b^2]\, \Phi+   [(n-c) b^3]\,+\nonumber\\
 &\phantom{2cm}& \Big([(n-c) a^2 +c^2]\,\Phi^2+2\,[(n-c)a b +cd]\,\Phi + [(n-c) b^2 +d^2]\Big)\,\Phi',\label{PhiR}
 \end{eqnarray}
 where we denote   $[.]= [.]\big|_{x=\Phi(t)} $.

 The system has the first order with respect to $e$ and the second order with respect to $\Phi$,
 it is supplemented by the initial condition for $e$,
 \begin{equation*}
 e(0)=-[E]^0
 \end{equation*}
and  boundary conditions for $\Phi$,
\begin{equation}\label{bc1}
\Phi(0)=\Phi_0,\quad \Phi(T_*)=x_+(T_*)=x_-(T_*).
\end{equation}

Now we
consider a situation where \eqref{cmcp} is not valid  for $t\in (t_1^*, t_2^*)$, $0< t_1^*< t_2^*\le T_*$, and at the points $t_1^*$ and $t_2^*$ condition \eqref{cmcp} becomes an equality, such that  $e(t_1^*)=e(t_2^*)=0$.
Then for $t\in(t_1^*, t_2^*) $  the rarefaction wave adjoint to the shock from the right side does not exist. Therefore from the right $V=V_+$ and $E=E_+$. Note that necessarily $\Phi(t_1^*)=x_\pm(t_1^*)$ and $e(t_1^*)=0$.  Assume  for definiteness that $\Phi(t_1^*)=x_+(t_1^*)$.

Thus, from \eqref{RH1}, \eqref{RH2} we obtain that the system for finding  $e(t)$ and $\Phi(t)$ looks like
 \begin{eqnarray}
 &e'\,&=\,([n]+c_-)\,\Phi' - (V_+ n_+ - ((n_--c_-)a_-\,\Phi+(n_--c_-)b_-),\nonumber
  \\
 & (e (\Phi')^2)'\,&=\,-n_+V_+^3+(n_--c_-)(a_-\Phi+b_-)^3\label{PhiR1}\\
 &&+(n_+V_+^2+E_+^2-(n_--c_-)(a_-\Phi+b_-)^2-(c_-\Phi+d_-)^2)\Phi',\nonumber
 %&& -   (n_--c_-) a_-^3\, \Phi^3 + 3 (n_--c_-) a_-^2 b_-\, \Phi^2 - 3 (n_--c_-) a_- b_-^2]\, \Phi+   (n_--c_-) b_-^3\,+\nonumber\\
 %&& \Big((n_--c_-) a_-^2 +c_-^2]\,\Phi^2+2\,((n_--c_-)a_- b_- +c_-d_-)\,\Phi + ((n_--c_-) b_-^2 +d_-^2)\Big)\,\Phi',\nonumber
 \end{eqnarray}
the boundary conditions are
\begin{equation*}
 e(t_1^*)=0,\quad \Phi(t_1^*)=x_+(t_1^*),\quad  \Phi(t_2^*)=x_+(t_2^*).
 \end{equation*}

In this case the boundary conditions for the step  $t\in (0, t_1^*)$ change from \eqref{bc1} to
\begin{equation*}\label{bc2}
\Phi(0)=\Phi_0,\quad \Phi(t_1^*)=x_+(t_1^*).
\end{equation*}
Condition \eqref{cmcp} now changes to
\begin{equation}\label{cmcp1}
c_-(t) \Phi(t)+ d_-(t) \ge E_+(t), \qquad t\in (t_1^*,t_2^*).
\end{equation}

If $\Phi(t_1^*)=x_-(t_1^*)$, then
 then system \eqref{RH1},~\eqref{RH2}
takes the form
\begin{eqnarray}
 &e'\,&=\,([n]-c_+)\,\Phi' + (V_- n_- - ((n_+-c_+)a_+\,\Phi+(n_+-c_+)b_+),\nonumber
  \\
 & (e (\Phi')^2)'\,&=\,n_-V_-^3-(n_+-c_+)(a_+\Phi+b_+)^3\label{PhiR2}\\
 &&-(n_-V_-^2+E_-^2-(n_+-c_+)(a_+\Phi+b_+)^2-(c_+\Phi+d_+)^2)\Phi',\nonumber
 %&& -   (n_--c_-) a_-^3\, \Phi^3 + 3 (n_--c_-) a_-^2 b_-\, \Phi^2 - 3 (n_--c_-) a_- b_-^2]\, \Phi+   (n_--c_-) b_-^3\,+\nonumber\\
 %&& \Big((n_--c_-) a_-^2 +c_-^2]\,\Phi^2+2\,((n_--c_-)a_- b_- +c_-d_-)\,\Phi + ((n_--c_-) b_-^2 +d_-^2)\Big)\,\Phi',\nonumber
 \end{eqnarray}
 the boundary conditions are
\begin{equation*}
 e(t_1^*)=0,\quad \Phi(t_1^*)=x_-(t_1^*),\quad  \Phi(t_2^*)=x_-(t_2^*).
 \end{equation*}

In this case the boundary conditions for the step  $t\in (0, t_1^*)$ change from \eqref{bc1} to
\begin{equation*}\label{bc2}
\Phi(0)=\Phi_0,\quad \Phi(t_1^*)=x_-(t_1^*).
\end{equation*}
and instead of condition \eqref{cmcp} we have
\begin{equation}\label{cmcp2}
c_+(t) \Phi(t)+ d_+(t) \le E_-(t), \qquad t\in (t_1^*,t_2^*).
\end{equation}

Note that condition \eqref{cmcp} can become an equality more than two time on $(0,T_*)$, and therefore a switching of rarefaction regions of different types takes place in every point of this kind (see Example 3 in Sec.\ref{S6}).

\subsection{Search for switching points}

The conditions \eqref{cmcp}, \eqref{cmcp1}, \eqref{cmcp2} help to define the boundaries of those regions on the characteristic plane $(t,x)$ in which $(e,\Phi)$ can be  found from \eqref{PhiR}, \eqref{PhiR1}, \eqref{PhiR2}, respectively.

Let us denote $\Psi_1(t)=-\frac{d_+(t)-d_-(t)}{c_+(t)-c_-(t)}$ the curve, obtained from \eqref{cmcp} with an equality and
assume that  $\Psi_1(t)$ changes sign only twice in the rarefaction region $(0, T_*)$, in points $t_1^*$ and $t_2^*$.
Note that $\Psi_1(t)$ can have more than two roots, see Sec.\ref{S4}.

 At the point of intersection of $\Phi(t)$ with $\Psi_1(t)$ the amplitude of the delta-function $e(t)$ vanishes and the singular shock between two adjoint rarefaction regions fails to exist. However, this point belongs to the characteristics $x_+(t)$ (or $x_-(t)$) and  the curve
 $\Psi_2^+(t)=\frac{E_+(t)-d_-(t)}{c_-(t)}$, defined by  \eqref{cmcp1} (or $\Psi_2^-(t)=\frac{E_-(t)-d_+(t)}{c_+(t)}$, defined by  \eqref{cmcp2}) passes through the same point. Therefore $t_1^*$ and $t_2^*$, where the singular shock between two adjoint rarefaction regions changes to the singular shock between one adjoint rarefaction regions and the state that does not depend of $x$ and vice versa are roots of the equation $\Psi_1(t)=\Psi_2^\pm(t)=x_\pm(t)$ and $e(t_1^*)=e(t_2^*)=0$. Thus, the points of switching between different types of singular shocks are $(t_1^*, x_\pm(t_1^*))$ and $(t_2^*, x_\pm(t_2^*))$ and we know in advance, where \eqref{PhiR}, \eqref{PhiR1} or \eqref{PhiR2} are valid.

\section{Conjugation of the shock waves in the compression and rarefaction regions}\label{S5}

Assume that the point $t=T_0$, $t=T_+$ and $t=T_-$  of the intersection of characteristics $x_+(t)$ and $x_-(t)$ is such that for $t\in (T_0, T_+)  $ we have shock wave and for $t\in (T_-, T_0)  $ we have rarefaction region, see \eqref{Kx}. Let $T_0=0$ and $x_\pm(0)=0$, we denote as $x=\Phi_\pm (t)$ the lines of discontinuity (singular shocks) for $t>0$ and $t<0$, respectively.

To construct shock we know $\Phi(0)=0$ and $\Phi(T_+)=x_\pm(T_+)$, so we have to solve the boundary problem for the second order equation
\eqref{Phidd}.

Let us show that if we {\it assume} that $\Phi(t)$ is $C^1$ -- smooth at zero, then
 we can find $\Phi'(0)$. If the rarefaction domain consists of two rarefaction waves, i.e. condition \eqref{cmcp} holds, then from \eqref{esw} and \eqref{eR}
\begin{eqnarray*}
&&e'(t)\\
&&=\,-[n V(t)]\big|_{x=\Phi_+(t)}+[n]\big|_{x=\Phi_+(t)} \Phi'(t) \\
&&=\,[n-c(t)]\big|_{x=\Phi_-(t)}\,\Phi'(t) -[(n-c(t))(a(t)\,\Phi(t)+b(t)]\big|_{x=\Phi_-(t)},
\end{eqnarray*}
therefore
\begin{eqnarray}
&&\Phi'(0)=\lim\limits_{t\to 0} \,\frac{[n V(t)]\big|_{x=\Phi_+(t)}-[(n-c(t)) (a(t)\,\Phi(t)+b(t)]\big|_{x=\Phi_-(t)}}{[c(t)]\big|_{x=\Phi_-(t)}}  \nonumber\\
&&= V^0_-\,+ C\,[V]^0, \label{Phidot}\\&& C=
\lim\limits_{t\to 0} \,\frac{c_+(t)}{c_+(t)-c_-(t)}= \frac{1}{1-\frac{B_+}{B_-}r},\quad r=\frac{\sqrt{n_-}}{\sqrt{n_+}}.\nonumber
\end{eqnarray}
%note that both numerator and denominator tend to infinity as $t\to 0$.

Similarly, if  condition \eqref{cmcp1} (or \eqref{cmcp2}) is satisfied for $t\in (T_-, T_0) $, then we obtain from
\eqref{esw} and \eqref{PhiR1} (or \eqref{PhiR2})
\begin{equation*}
\Phi'(0)= V^0_-
\end{equation*}
or
\begin{equation*}
\Phi'(0)= V^0_+,
\end{equation*}
respectively.

Since for the shock
$e(t)$ is known from \eqref{ee}, we can find
$$
e(T_+)=-[E(T_+)]+[n]\Phi(T_+)=e(0).
$$

Now {\it assume} that the problem is periodical, $L=T_++T_-$ is a period, and there is only one switch from a shock to rarefaction with two rarefaction waves (and vice versa) within the period.

Then $\Phi(T_+)=\Phi(T_-)$, $e(T_+)=e(T_-)$ and
$$
\Phi'(T_+)=\Phi'(T_-)= V_-(T_\pm)\,+ C\,[V(T_\pm)], \qquad C=
\lim\limits_{t\to T_\pm} \,\frac{c_+(t)}{c_+(t)-c_-(t)}.
$$

Thus, we can solve equation \eqref{Phidd} in the case of a shock wave and system \eqref{eR}, \eqref{PhiR} in the case of a rarefaction subject  to initial conditions,
which is easier numerically.

Let us show, however, that in the general case one cannot require the smoothness of the singular shock wave at the switching point from compression to rarefaction. Indeed, that condition \eqref{accept} for the shock starting from zero to the right  dictates $V_+^0\le \Phi'(0)\le V_-^0$, therefore
\eqref{Phidot} implies
\begin{eqnarray}
&&V_+^0\le V^0_-\,+ C\,[V]^0\le V_-^0,\qquad 0\le C\le 1, \qquad \frac{B_+}{B_-}\le 0,\nonumber \\
&&\frac{\tan{\frac{\sqrt{n_+}T_*}{2}}}{\tan{\frac{\sqrt{n_-}T_*}{2}}}\le 0. \label{Phin}
\end{eqnarray}
Evidently, that if ${n_-}\to {n_+}$, then \eqref{Phin} cannot be satisfied, therefore the curve $x=\Phi(t)$ should lose smoothness at the point of conjugation of singular shocks inside and outside the rarefaction regions.

The conditions for the conjugation of shock waves inside a rarefaction region can be obtained in a completely analogous way from  \eqref{PhiR} and \eqref{PhiR1} (or \eqref{PhiR2}). Under the condition of smoothness of $\Phi(t)$, when transitioning from a rarefaction region with two rarefaction waves to a region with one rarefaction wave at a point $t_*$, the condition 
\begin{equation*}
\Phi'(t_*)= V_\pm(t_*)
\end{equation*}
must be satisfied  if the transition occurs on the characteristic $x_\pm(t)$.

Note that even if the problem is periodic, there may be several switching points between the shock wave and the rarefaction region within a period. Then, to find the derivative $\Phi'$ at an arbitrary intersection point $(t_*,x_*)$ of the characteristics $x=x_\pm (t)$, we must shift it to the origin and apply the formulas obtained above with the  change of  variables $t-t_*$ and $x-x_*$.

However, both equation \eqref{Phidd}  and system \eqref{eR}, \eqref{PhiR} degenerate at the points where $\Phi'(t)=0$, therefore the numerical solution is challenging if $\Phi(t)$ loose monotonicity. Nevertheless, $\Phi(t)$ is basically not monotonic.
Therefore, constructing singular discontinuity lines is, in any case, a very difficult computational problem.

\section{Examples}\label{S6}

In this section, we  provide several examples to illustrate our reasoning. On the plane of characteristics $(t,x)$ the red line denotes a singular shock wave, the interface between media with different base densities. In Fig.1 on the left, we show the arrangement of the characteristics $x_+(t)$ and $x_-(t)$ for a case where the background densities are equal. The singular shock wave and the rarefaction wave alternate periodically, and there is no dividing line within the rarefaction wave.
Fig.1 on the right, shows an example of the arrangement of the characteristics $x_+(t)$ and $x_-(t)$ for a case of incommensurate oscillation periods. We see that compression and rarefaction zones alternate in a complex manner. Rarefaction zones are indicated by solid fill.  Within the rarefaction zone, the interface can be adjoined by either two rarefaction waves on the left and right, or by a rarefaction wave on one side and a state independent of $x$ on the other. The position of the interface is shown schematically.

\begin{center}
\begin{figure}[htb!]
\hspace{-1cm}
%\begin{minipage}{0.3\columnwidth}
%\centerline{
%\includegraphics[scale=0.25]{Ex2V3.ps}
%\vspace{-0.5 cm}
%\end{minipage}
\hspace{1cm}
\begin{minipage}{0.4\columnwidth}
%\centerline{
\includegraphics[scale=0.35]{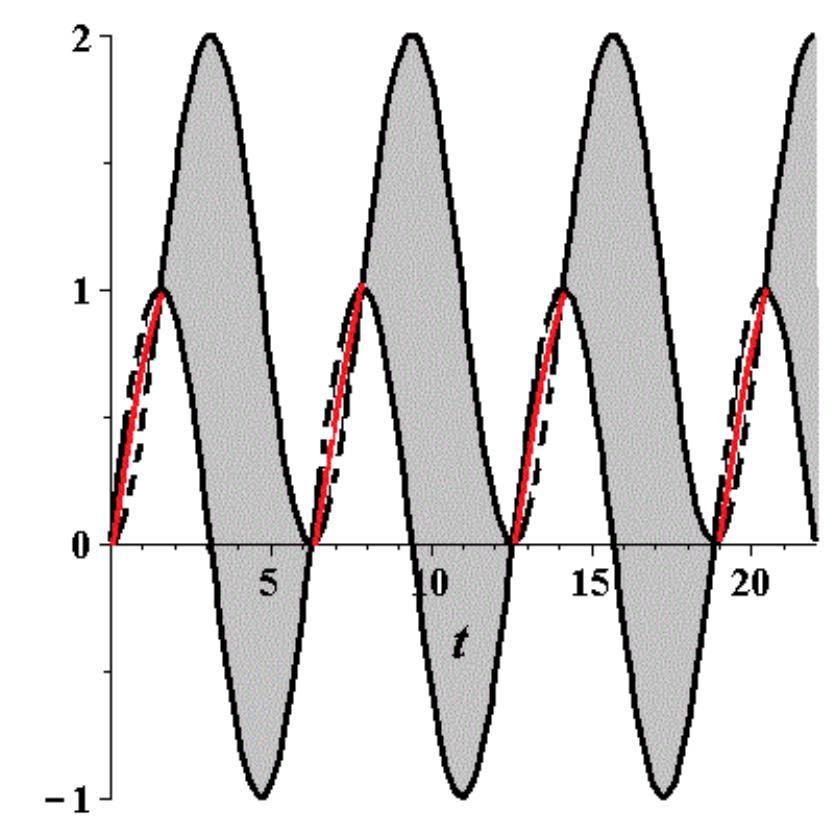}
%\vspace{-0.5 cm}
%\caption{ }%\label{Pic2}
\end{minipage}
\hspace{1.5cm}
\begin{minipage}{0.4\columnwidth}
%\centerline{
\includegraphics[scale=0.35]{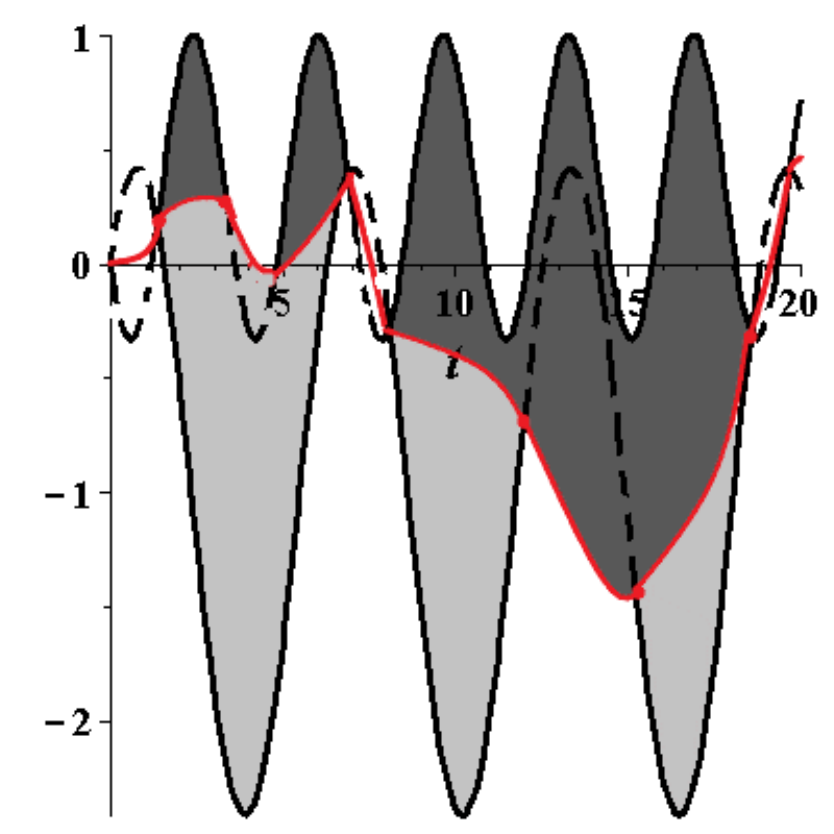}
%\vspace{-0.5 cm}
%\caption{$n_-=3, n_+=1,  V_+^0=-1,  V_-^0=1,   E_+^0=-1,  E_-^0=1.$}\label{Pic2}
\end{minipage}
\caption{Plane $(t,x)$. Characteristics $x_-(t), x_+(t)$ and the position of the interface $\Phi(t)$. Left: Example 1, $n_-=n_+=1,  V_+^0=0,  V_-^0=1,   E_+^0=-1,  E_-^0=0.$  Right: Example 2, $n_-=1, n_+=3,  V_+^0=-1,  V_-^0=1,   E_+^0=-1,  E_-^0=1.$ }\label{Pic2}
\end{figure}
\end{center}

Fig.2  shows  the situation with different background densities in the case of commensurate oscillation periods and, accordingly, periodic alternation of compression and rarefaction regions (Example 3). The intersection of the characteristics $x_-$ and $x_+$ are in the points $T_0=2 k \pi$ and $T_*=1.035895953 + 2 k \pi$, $k\in \mathbb Z$. %The condition \eqref{Tstar} is not met.
% $B_+\sin{\sqrt{n_+}t}+\cos{\sqrt{n_+}t}-1=0  $
%The dotted line represents the graph of the denominator $c_+(t)$.
The intersection of $x_+(t)$ with $\Psi_2^+(t)$ within the rarefaction region $t\in (T_*, 2 \pi)$ are in the points $t_1^*=2.176190164 + 2 k \pi$,  $t_2^*=3.920405792 + 2 k \pi$, and $t_3^*=5.916224372 + 2 k \pi$, $k\in \mathbb Z$. In these points the zone where two rarefaction waves adjoin the boundary of the media $\Phi$ switches to a zone where only one such wave occurs, and vice versa. The position of the interface  $x=\Phi(t)$  is indicated schematically (Fig.2, left). Fig.2, right, shows the origin of the switching points as intersection of characteristic $x_+$ and curves $\Psi_1$ (dash) and $\Psi_2^+$ (dot), boundaries of the domains where conditions \eqref{cmcp} and \eqref{cmcp1} hold.

Figs.3 and 4 present a detailed picture of position of the singular interface $\Phi$ and the structure of the rarefaction region.

\begin{center}
\begin{figure}[htb!]
\hspace{-1cm}
%\begin{minipage}{0.3\columnwidth}
%\centerline{
%\includegraphics[scale=0.25]{Ex2V3.ps}
%\vspace{-0.5 cm}
%\end{minipage}
\hspace{1cm}
\begin{minipage}{0.4\columnwidth}
%\centerline{
\includegraphics[scale=0.35]{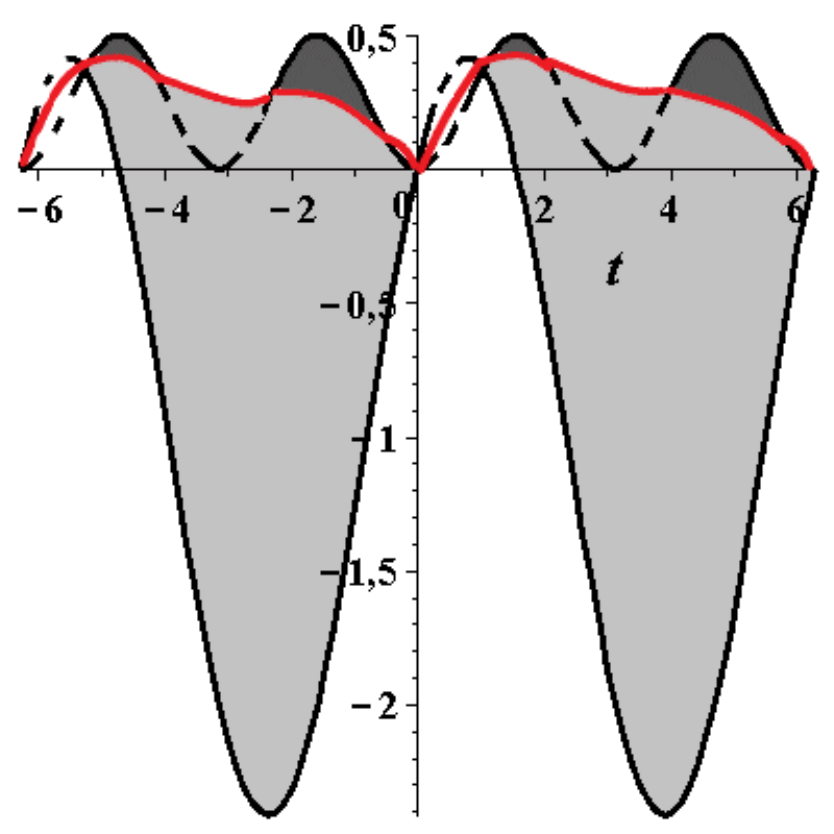}
%\vspace{-0.5 cm}
%\caption{ }%\label{Pic2}
\end{minipage}
\hspace{1.5cm}
\begin{minipage}{0.4\columnwidth}
%\centerline{
\includegraphics[scale=0.35]{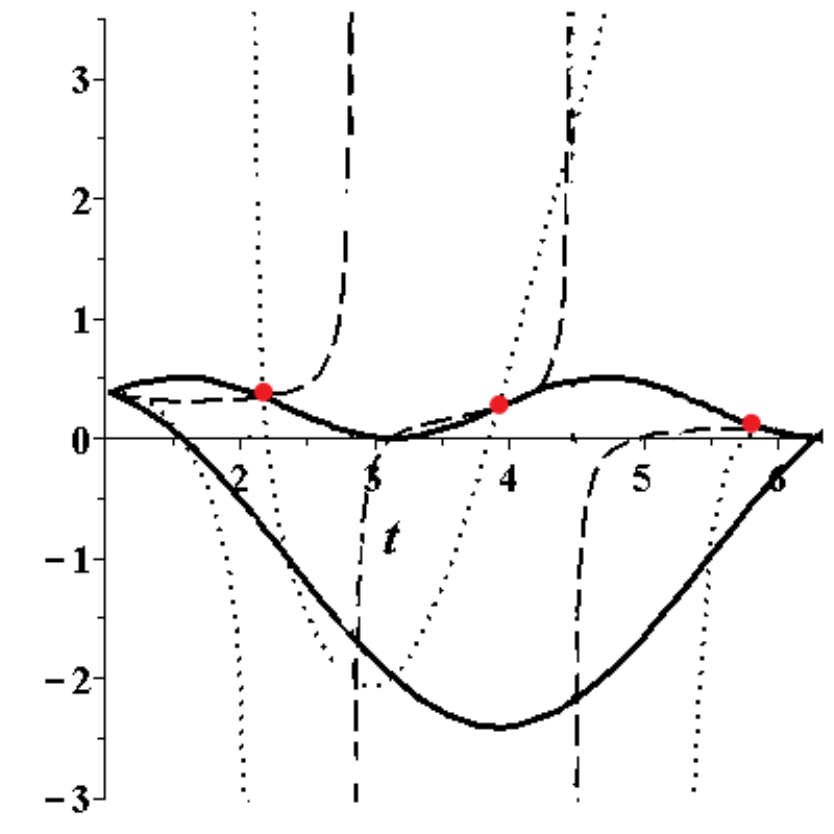}
%\vspace{-0.5 cm}
%\caption{$n_-=3, n_+=1,  V_+^0=-1,  V_-^0=1,   E_+^0=-1,  E_-^0=1.$}\label{Pic2}
\end{minipage}
\caption{Plane  $(t,x)$, Example 3: $n_-=1, n_+=4,  V_+^0=0,  V_-^0=1,   E_+^0=-1,  E_-^0=1.$ Left:  a schematic position of the interface. Right: rarefaction region within one period, graphs of $\Psi_1(t)$ (dash) and $\Psi_2^+(t)$ (dot), switching points $t_1^*$, $t_2^*$, $t_3^*$ }\label{Pic2}
\end{figure}
\end{center}

\begin{center}
\begin{figure}[htb!]
\hspace{-1cm}
%\begin{minipage}{0.3\columnwidth}
%\centerline{
%\includegraphics[scale=0.25]{Ex2V3.ps}
%\vspace{-0.5 cm}
%\end{minipage}
\hspace{1cm}
\begin{minipage}{0.4\columnwidth}
%\centerline{
\includegraphics[scale=0.35]{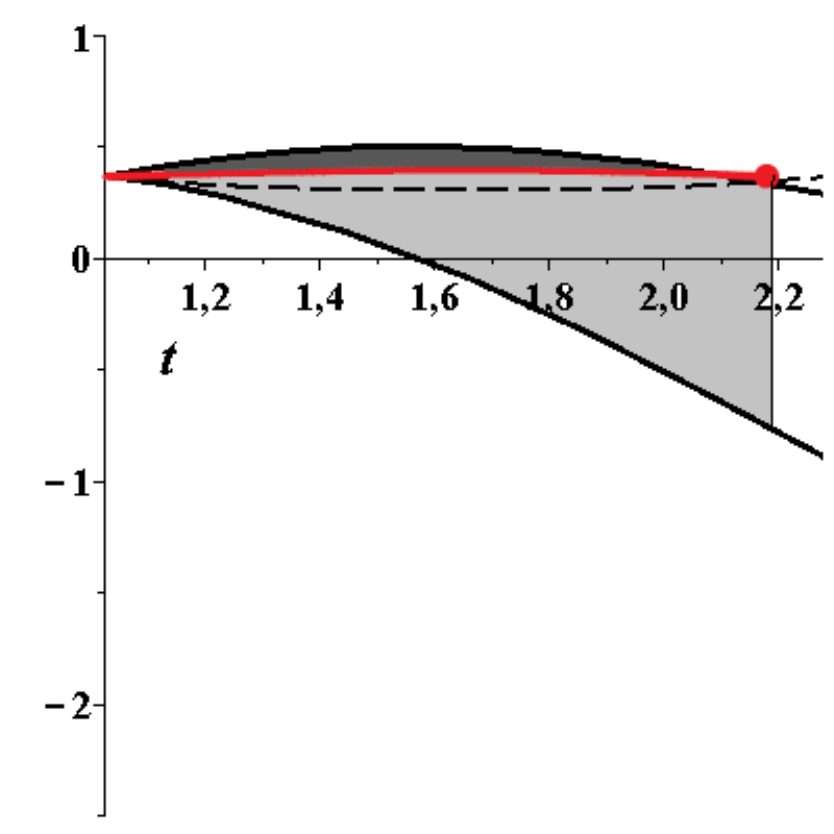}
%\vspace{-0.5 cm}
%\caption{ }%\label{Pic2}
\end{minipage}
\hspace{1.5cm}
\begin{minipage}{0.4\columnwidth}
%\centerline{
\includegraphics[scale=0.35]{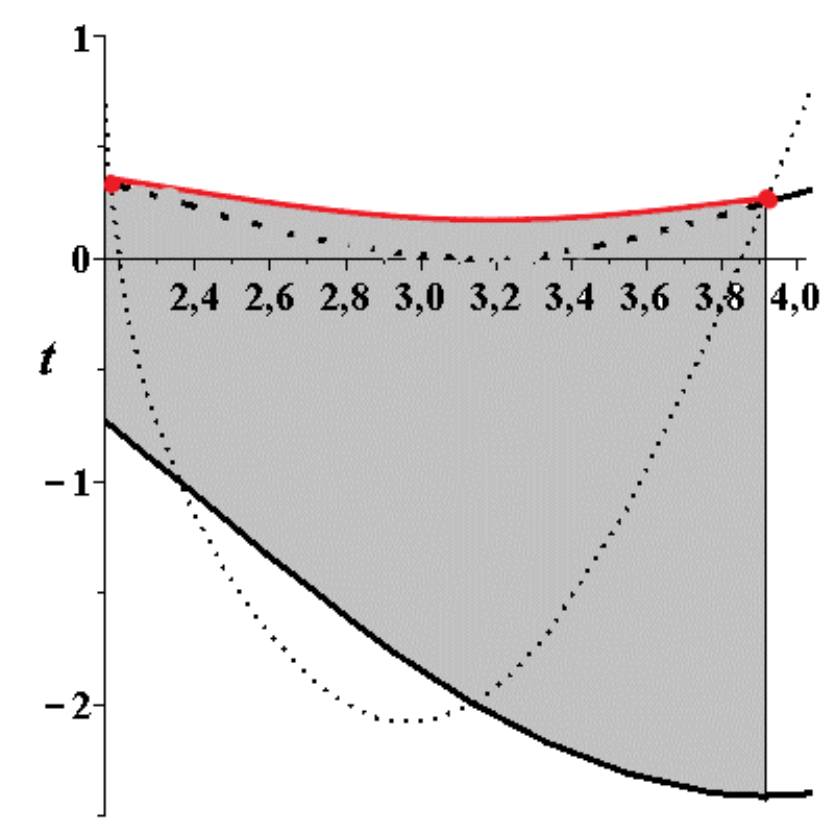}
%\vspace{-0.5 cm}
%\caption{$n_-=3, n_+=1,  V_+^0=-1,  V_-^0=1,   E_+^0=-1,  E_-^0=1.$}\label{Pic2}
\end{minipage}
\caption{Plane $(t,x)$, Example 3, detailed picture of
characteristics $x_-(t), x_+(t)$ and the position of the interface $\Phi(t)$.
  Left:  $t\in (0, t_1^*)$, rarefaction region with two rarefaction waves. Right: $t\in (t_1^*, t_2^*)$, rarefaction region with one rarefaction wave; the characteristic $x_+$ is removed.}\label{Pic3}
\end{figure}
\end{center}
\begin{center}
\begin{figure}[htb!]
\hspace{-1cm}
%\begin{minipage}{0.3\columnwidth}
%\centerline{
%\includegraphics[scale=0.25]{Ex2V3.ps}
%\vspace{-0.5 cm}
%\end{minipage}
\hspace{1cm}
\begin{minipage}{0.4\columnwidth}
%\centerline{
\includegraphics[scale=0.35]{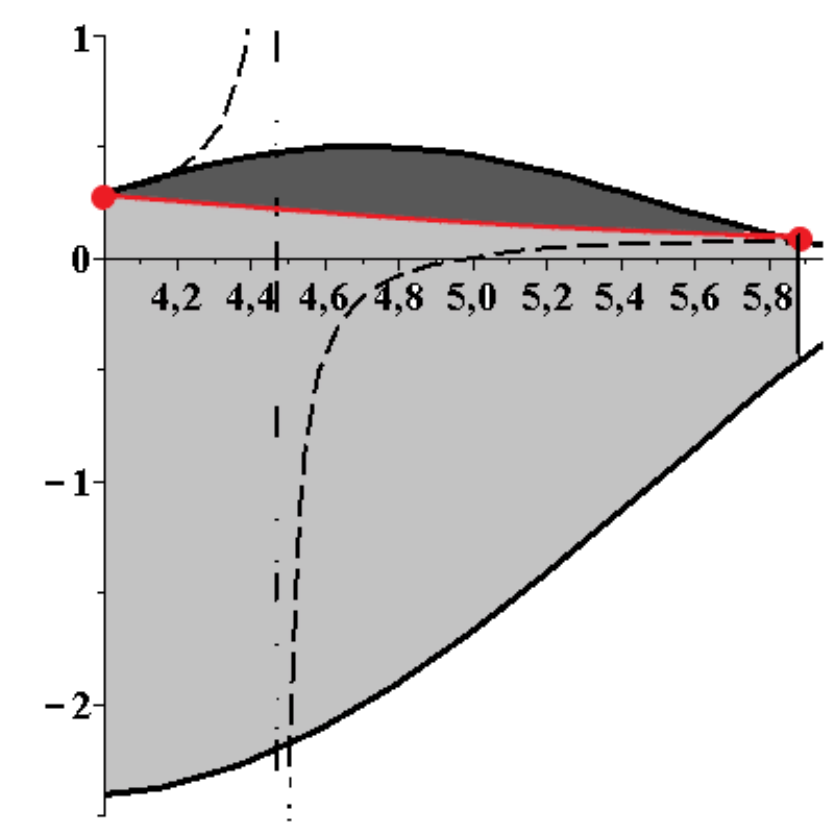}
%\vspace{-0.5 cm}
%\caption{ }%\label{Pic2}
\end{minipage}
\hspace{1.5cm}
\begin{minipage}{0.4\columnwidth}
%\centerline{
\includegraphics[scale=0.35]{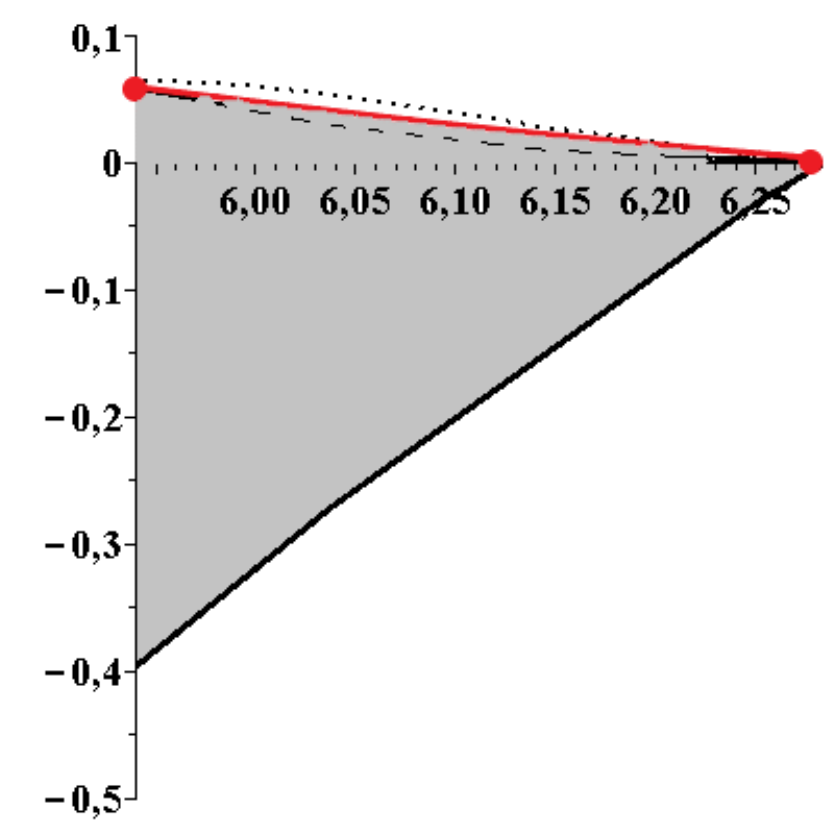}
%\vspace{-0.5 cm}
%\caption{$n_-=3, n_+=1,  V_+^0=-1,  V_-^0=1,   E_+^0=-1,  E_-^0=1.$}\label{Pic2}
\end{minipage}
\caption{Plane $(t,x)$,  Example 3,  detailed picture of characteristics $x_-(t), x_+(t)$ and the position of the interface $\Phi(t)$.  Left:  $t\in (t_2^*, t_3^*)$, rarefaction region with two rarefaction waves. Right: $t\in (t_3^*, 2 \pi)$, rarefaction region with one rarefaction wave; the characteristic $x_+$ is removed.  }\label{Pic4}
\end{figure}
\end{center}

Note that the position of $x=\Phi(t)$ for $t=[0, T_*]$ can be easily found numerically by the shooting method, since in this case $\Phi'(t)\ne 0$. The computation of the position of $x=\Phi(t)$ for $t=[T_*, t_1^*]$, $t=[t_1^*, t_2^*]$, $t=[t_2^*, t_3^*]$ and $t=[t_3^*, 2\pi]$
are significantly complicated by the presence of degeneracy points in these regions, where $\Phi'(t)=0$ and $e(t)=0$. Note that theoretical results on the existence and uniqueness of such boundary value problems are unknown.

\section{ Conclusion and discussion}\label{S7}

We  considered a one-dimensional case of cold plasma equations involving a junction of two media with constant but different background densities. The interface between the media is a strongly singular discontinuity, where the generalized Rankine-Hugoniot conditions must be satisfied. The problem is to find the values of the solution components (velocity, density, and electric field strength) to the left and right of the discontinuity, as well as the position of the discontinuity itself corresponding to the interface between the two media. Thus, we solve the Riemann problem. The initial data are assumed to be constant to the left and right of the discontinuity.

When there is no separation of the media, the Riemann problem is solved in \cite{Riemann_Roz}. In this case, its solution consists of alternating singular shock waves and continuous rarefaction waves. It can be shown that the solution always exists and, under certain conditions, is unique.  In \cite{GRT} the Riemann problem was considered for all variants of the Euler-Poisson equations, which yield cold plasma equations in the particular case of a repulsive interaction force between particles and a nonzero background density. In all cases, the solution has interesting distinctive features, but is nevertheless constructed in a fairly standard manner.

The assumption of a discontinuity in the background density significantly complicates the problem. It turns out that the rarefaction region must also contain a singular shock wave. Depending on the ratio $r=\frac{\sqrt{n_-}}{\sqrt{n_+}}$, where $n_\pm$ are the background density values to the right and left of the discontinuity, the compression and rarefaction regions can alternate intricately.
In particular, if $r$ is rational, meaning the oscillations of the media to the left and right of the discontinuity are commensurate, then the problem becomes periodic. However, within this period, there can be quite a few switches from the singular shock wave to the rarefaction region and vice versa.
However, the rarefaction region itself can have a variety of structures. Specifically, a singular interface can be approached by two waves on each side, or by a rarefaction wave on only one side. Within each rarefaction region, such structures alternate.

Although in this paper we have posed problems that formally determine the solution, many unanswered questions remain. In particular, the question of the existence and uniqueness of a solution remains open in the general case. It may well turn out that the singular interface between the media cannot be smooth in some cases, and the uniqueness of the solution is open in the case if the interface loses its monotonicity. Furthermore, finding a solution numerically becomes a real challenge due to the presence of degeneracy points.

Thus, a seemingly insignificant modification to a fairly standard problem has led to the emergence of a whole new field of research.

\section{Acknowledgments}
Supported by RSF grant 23-11-00056 through RUDN University.

\end{document}